\documentclass[10pt]{amsart}

\usepackage[english]{babel}

\usepackage[utf8]{inputenc}

\usepackage{amsmath,amsthm,verbatim,amssymb,amsfonts,amscd, graphicx}
\usepackage{graphics}
\usepackage{tensor}
\usepackage{enumitem}
\usepackage[utf8]{inputenc}
\usepackage{graphicx}
\usepackage{xcolor}
\usepackage{tikz-cd}
\usepackage{mathtools}
\usepackage{esint}

\usepackage[english]{babel}
\usepackage{graphicx}
\usepackage{subcaption}
\usepackage[T1]{fontenc}
\usepackage{hyperref}
\usepackage{amsmath}
\usepackage{amssymb}
\usepackage{amsfonts}
\usepackage{amssymb}
\usepackage{amsthm}
\usepackage{amscd}
\usepackage{graphicx}
\usepackage{xcolor}
\usepackage{float}
\usepackage{slashed}
\usepackage{tikz}

\usepackage{geometry}

\theoremstyle{plain}
\newtheorem{theorem}{Theorem}[section]

\newtheorem{lemma}[theorem]{Lemma}
\newtheorem{remark}[theorem]{Remark}

\theoremstyle{definition}
\newtheorem{definition}[theorem]{Definition}
\theoremstyle{definition}
\newtheorem{example}[theorem]{Example}

\def\R{\mathbb{R}}

\def\inf{\operatorname{inf}}

\title{The Riemannian Penrose Inequality with Matter Density}

\author{Hubert Bray}
\address{Department of Mathematics, Duke University, Durham, NC, 27708, USA}
\email{bray@math.duke.edu}
\author{Yiyue Zhang}
\address{Department of Mathematics, University of California, Irvine, CA, 92697, USA}
\email{yiyuez4@uci.edu}

\begin{document}

\maketitle

\begin{abstract}
    Riemannian Penrose Inequalities are precise geometric statements that imply that the total mass of a zero second fundamental form slice of a spacetime is at least the mass contributed by the black holes, assuming that the spacetime has nonnegative matter density everywhere. In this paper, we remove this last assumption,  
    and prove stronger statements that the total mass is at least the mass contributed by the black holes, plus a contribution coming from the matter density along the slice. 
    
    We use the first author's conformal flow to achieve this, combined with Stern's harmonic level set techniques in the first case, and spinors in the second case. We then compare these new results to results previously known from Huisken-Ilmanen's inverse mean curvature flow techniques.
\end{abstract}

\section{Introduction}
The Penrose conjecture is one of the most important open problems in mathematical relativity.
A comprehensive survey about the Penrose conjecture for spacelike slices of spacetimes can be found in \cite{Penrose Survey}.  When the slice is assumed to have zero second fundamental form in the spacetime, the Penrose conjecture is known as the Riemannian Penrose inequality, 
resolved by two different approaches: inverse mean curvature flow (IMCF) \cite{HI01} which works for a single black hole, and the first author's conformal flow \cite{bray01} which works for any number of black holes, and in space-like dimensions up to 7 \cite{bray2009riemannian}.  The latter approach
is also used to prove 
the charged Riemannian Penrose inequality in \cite{charge Penrose}.

In dimension 3, which will be the focus of this paper, the Riemannian Penrose inequality asserts that the total mass $m$ is greater than $\sqrt{A/16\pi}$, where $A$ is the area of the outermost minimal surface of the slice $M^3$. In this paper, we will prove stronger inequalities of the form
\begin{equation} \label{1}
    m\ge \sqrt{\frac{A}{16\pi}}+\frac{1}{16\pi}\int_{M^3} \left(RQ(x)+P(x)\right)dV,
\end{equation}
where $R$ is the scalar curvature of $M^3$ and $Q(x)$ and $P(x)$ are nonnegative functions which vanish inside the outermost minimal surface (representing the black holes). Hence, $R\ge 0$ implies the usual Riemannian Penrose inequality. Also note that $R/{16\pi}$ is matter density, so integrating it is quite natural, especially since the $Q(x)$ we discuss all go to one at infinity. The fact that $Q(x)$ is typically less than one represents a negative potential energy contribution, which is also expected. Understanding the physical interpretation of these inequalities is a good reason to study them.

An even more exciting reason is this: Every inequality of the form of (\ref{1}) defines a system of PDE's which, when solutions with certain properties exist, imply the Penrose Conjecture \cite[page 40]{PDE Penrose}. Hence, every inequality of the form of (\ref{1}) gives a plan of attack for the Penrose conjecture.


The easiest way to prove an inequality of the form of (\ref{1}) is to use inverse mean curvature flow \cite{HI01}, as demonstrated in  \cite[Theorem 7]{PDE Penrose}. However, there are two main limitations with the resulting formula. First, the inequality is only proven when the outermost minimal surface is connected, or if $A$ is defined to be the area of the largest connected component of the outermost minimal surface. Second, while $Q(x) \ge 0$, $Q(x) = 0$ in the jump regions of inverse mean curvature flow, which destroys the ellipticity of the system of PDE's defined above.
Hence, while the achieved formulas are very nice, there are good reasons to look for more formulas of this form.

In this paper, we will present two more inequalities of the form of (\ref{1}) above using the conformal flow. 
A key advantage of these inequalities arising from the conformal flow is that they guarantee $Q(x)$ to be positive almost everywhere outside the outermost minimal surface. Additionally, the conformal flow approach works for any number of black holes.

For convinence, we only study harmonically flat manifolds which is a representative case of asymptotically flat manifolds introduced in \cite{bray01}, since  every asymtotically flat metric is a harmonically flat metric up to a $C^0$ perturbation.

\begin{definition}
An smooth complete manifold $M^n$ is \emph{harmonically flat}  if there exists a compact set $K$ such that $M^n\setminus K$ is a disjoint union of multiple ends and each end $E_k$ is diffeomorphic to $\R^n\setminus B_1(0)$. Moreover,  on each end $E_k$, there exists a harmonic function $u$ asymptotic to a constant such that 
\[g_{ij}=u^\frac{4}{n-2}\delta_{ij}.\]
\end{definition}

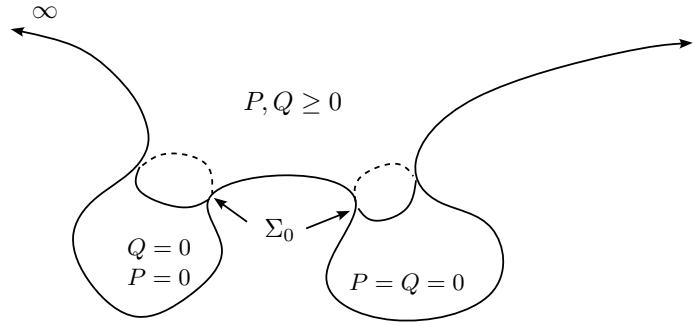
\begin{figure}[H]
    \centering
    \includegraphics[width=0.6\textwidth]{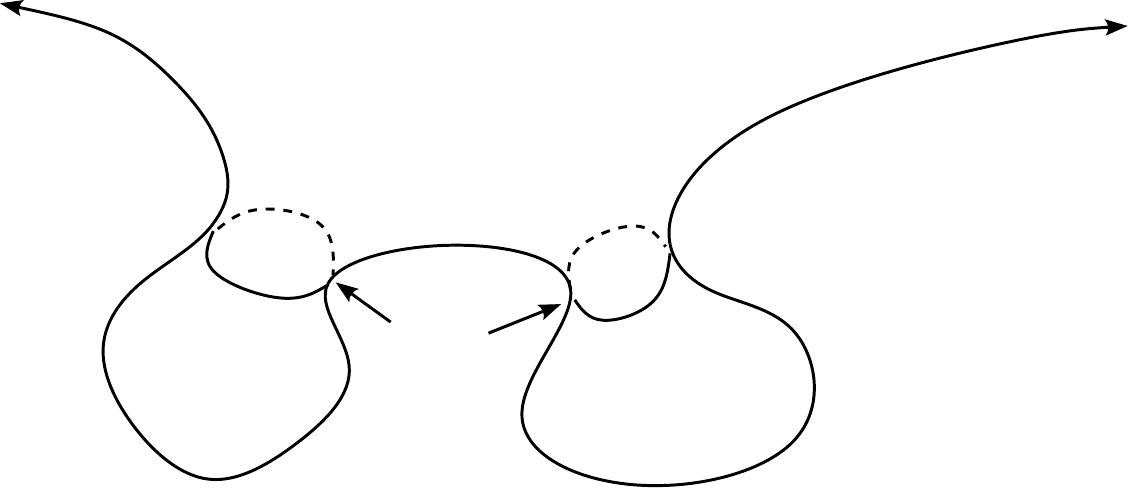}
    \put(-162,32){$\Sigma_0$}
    \put(-214,14){\small{$P=0$} }
     \put(-214,25){\small{$Q=0$} }
    \put(-250,115){$\infty$}
    \put(-170,80){$P,Q\ge0$}
    \put(-130,12){\small{$P=Q=0$}}
    \caption{An asymptotically flat manifold  with a disconnected outermost minimal hypersurface $\Sigma_0$}
    \label{fig:enter-label}
\end{figure}

Note that a harmonically flat end has zero scalar curvature and $u$ satisfies 
\[u(x)=a+b|x|^{2-n}+O(|x|^{1-n}),\]
where $a$ and $b$ are some constants. Then the ADM mass of the end $E_k$ is $2ab$.

Here is the definition of the conformal flow in dimension $3$. 
\begin{definition} \label{conflow}
    Let $(M^3,g)$ be a harmonically flat manifold. $(M^3,g_t,\Sigma(t))$ is a solution to the conformal flow if it satisfies
    \begin{itemize}
    \item $(M^3,g_t)$ is a smooth harmonically flat manifold outside $\Sigma_t$. 
    \item $\Sigma(t)$ is an outermost minimal surface in $(M^3,g_t)$.
        \item $g_t=u^4_tg$, where 
        \begin{equation} \label{ut}
    u_t(x)=1+\int_0^tv_s(x)dt
\end{equation}
and $v_t(x)$ satisfies
\begin{equation}
\begin{cases}
    \Delta_{g} v_t(x)\equiv 0\quad &\textit{outside }\Sigma(t) 
    \\ v_t(x)=0 \quad &\textit{on and inside } \Sigma(t) 
    \\ \lim_{x\to\infty} v_t(x)=-e^{-t}
\end{cases}    
\end{equation}
    \end{itemize} 
\end{definition}
Let $A$ be the area of the outermost minimal surface of $(M^3,g)$. 
Let $M^3_{\Sigma(t)}$ be the connected component of $M^3\setminus\Sigma(t)$ containing the chosen end  
and $t(x):=\inf\{t| x\notin M^3_{\Sigma(t)}\}$. By definition, $u_t$ is harmonic on $M^3_{\Sigma(t)}$. 

The theorem presented below establishes inequality \eqref{1} with new $P(x)$ and $Q(x)$. 
 This theorem is obtained by using the conformal flow \cite{bray01} and Stern's harmonic level set method \cite{BKKS,Stern}.
\begin{theorem}\label{harmonic theorem}
Let $(M^3,g)$ be a 3-dimensional smooth complete harmonically flat manifold. Define $m$ as the ADM mass of a chosen end and let $R$ be the scalar curvature of $(M^3,g)$. Suppose the outermost minimal surface of $(M^3,g)$ is nonempty and connected. Let $p_t$ and $q_t$ be two harmonic functions on $(M^3,g)$ satisfying
\begin{align}
     \Delta p_t=0,&\quad p_t\to e^{-3t}x_1 \text{ at } \infty \text{ and }  \nu_t(p_t)=p_t\nu_t(\log u_t) \text{ at }\Sigma(t)
     \\ \Delta q_t=0,&\quad q_t\to e^{-3t}x_1 \text{ at } \infty \text{ and }  q_t=0 \text{ at }\Sigma(t)
 \end{align} 
  where $u_t$ and $\Sigma(t)$ are defined in Definition \ref{conflow},  $\nu_t$ is the unit normal vector on $\Sigma(t)$ pointing towards to the chosen end with respect to $g$, $x_1$ is the first component of $x$ in the asymptotically flat coordinate on the chosen end. 

 Let $\phi_t$ be the harmonic function satisfying
 \begin{equation}
    \Delta \phi_t=0;\quad \lim_{x\to \infty}\phi_t=e^{-t};\quad \phi_t=\frac{u_t}{2} \textit{ on } \Sigma_t,  
\end{equation}
and denote $\mathcal{U}_{t,1}=\phi_t^{-1}(p_t+q_t)$, $\mathcal{U}_{t,2}=(u_t-\phi_t)^{-1}(p_t-q_t)$,  then 
        \begin{equation}
            \begin{split}
                Q(x):=&\int_0^{t(x)} (|\nabla\mathcal{U}_{t,1}|+|\nabla \mathcal{U}_{t,2}|)dt,
                \\ P(x):=&\sum_{j=1}^2\int_0^{t(x)}\frac{\left|\nabla^2 \mathcal{U}_{t,j}-2\phi_t^{-1}\nabla \phi_t\otimes\nabla \mathcal{U}_{t,j}-2\phi_t^{-1}\nabla \mathcal{U}_{t,j}\otimes\nabla \phi_t+2\phi_t^{-1}\langle\nabla\phi_t,\nabla \mathcal{U}_{t,j}\rangle g\right|^2}{|\nabla \mathcal{U}_{t,j}|},
            \end{split}
        \end{equation}
        and
       \begin{equation} 
    m\ge \sqrt{\frac{A}{16\pi}}+\frac{1}{16\pi}\int_{M^3} \left(RQ(x)+P(x)\right)dV.
\end{equation} 
\end{theorem}
\begin{remark}
    We need to assume the outermost minimal surface is connected so that we can apply the integral formula using harmonic level method. This assumption is necessary as the harmonic level method relies on the Gauss-Bonnet theorem.
\end{remark}
  In particular, we compute $Q(x)$ in Example \ref{CFS} under the Schwarzschild metric. Similarly, IMCF also gives the same answer \cite{PDE Penrose}. 
\begin{example} \label{imcf}
    Let $(M^3,g)$ be the spatial Schwarzschild manifold, i.e.,  $g_{ij}=(1+\frac{m}{2r})^4\delta_{ij}$.  Under
    IMCF \cite[Page 38]{PDE Penrose}, the mean curvature of the hypersurface $\Sigma_r:=\{|x|=r\}$ is $H=(1+\frac{m}{2r})^{-3}(\frac{2}{r}-\frac{m}{r^2})$, $|\Sigma_r|=4\pi r^2(1+\frac{m}{2r})^4$, then $Q=H\sqrt{\frac{|\Sigma_t|}{16\pi}}=(1+\frac{m}{2r})^{-1}(1-\frac{m}{2r})$. 
    
    Since the Schwarzschild metric is the rigidity case of the Riemannian Penrose inequality, we have $P(x)=0$.
\end{example}

Spinors can also be used to prove the Positive Mass theorem and provide a mass formula \cite{Witten81, PT81}.   A second version of $P(x)$ and $Q(x)$ are defined by spinors in the following theorem, which produces an equality in \eqref{2}.  A similar formula is established in \cite[Equation (3.2) and (3.3)]{HanKhuri}.
\begin{theorem} \label{spinor theorem}
    Let $(M^3,g)$ be a 3-dimensional complete harmonically flat manifold with nonnegative scalar curvature $R$. Let $t(x)$ be the minimal time such that $x$ does not belong to $M^3_{\Sigma(t(x))}$. 
    Let $\mathbf{p}_t$ and $\mathbf{q}_t$ be the harmonic spinors satsifying
    \begin{equation}
        \begin{split}
            &\mathcal{D}(\mathbf{p}_t)=0, \quad
           i\nu_t\cdot \mathbf{p}_t=\mathbf{p}_t \text{ on }\Sigma(t)
           \text{ and } 
            \lim_{x\to\infty} \mathbf{p}_t=e^{-2t}\psi_0 ;
            \\ &\mathcal{D}(\mathbf{q}_t)=0, \quad i\nu_t\cdot \mathbf{q}_t=-\mathbf{q}_t \text{ on }\Sigma(t)\text{ and } \lim_{x\to\infty} \mathbf{q}_t=e^{-2t}\psi_0,
        \end{split}
    \end{equation}
    where $\psi_0$ is the unit length constant spinor at $
    \infty$. 
    Denote $\psi_t^l:=\frac{1}{2}(\mathbf{p}_t+(-1)^l\mathbf{q}_t)$, $l=0,1$; and $\phi_{t,l}$ are harmonic functions satisfying
    \begin{equation}
        \Delta \phi_{t,l}=0, 
        \quad 
        \lim_{x\to\infty}\phi_{t,l}=\frac{1+(-1)^l}{2e^{t}} \text{ and }\phi_{t,l}=\frac{u_t}{2} \text{ on }\Sigma_t,
    \end{equation}
    where $u_t$ is defined in Equation \eqref{ut}.
    Then 
    \begin{equation} \label{PQ2}
        \begin{split}
            Q(x):=&2\int_0^{t(x)}
( \phi_{t,0}^{-2} |\psi^0_t|^2+\phi_{t,1}^{-2} |\psi^1_t|^2)dt,
\\ P(x):=&8\int_0^{t(x)}\sum_{l=0}^1\phi_{t,l}^2\sum_{j=1}^3\left|\nabla_{e_j}(\phi_{t,l}^{-2}\psi_t^l)-[\nabla_{e_j} \log \phi_{t,l}+{e_j}\cdot \nabla \log \phi_{t,l}](\phi_{t,l}^{-2}\psi_t^l)\right|^2dt,
        \end{split}
    \end{equation}
    and
    \begin{equation}\label{2}
         m= \sqrt{\frac{A}{16\pi}}+\frac{1}{16\pi}\int_{M^3} \left(RQ(x)+P(x)\right)dV.
    \end{equation}
\end{theorem}
\begin{remark}
    We presume the scalar curvature to be nonnegative to ensure the existence of harmonic spinors. 
    It remains uncertain whether a harmonic spinor exists on asymptotically flat manifolds without the assumption of nonnegative scalar curvature.
    However, for the purpose of proposing a proof for the Penrose Conjecture in \cite{PDE Penrose},  the scalar curvature could be potentially negative, but the existence theory for the harmonic spinors would have to be established, perhaps using the dominant energy condition somehow.  
\end{remark}

\section{The Conformal flow}

\begin{figure}[H]
\centering
\begin{subfigure}{0.3\textwidth}
\centering
\includegraphics[width = \textwidth]{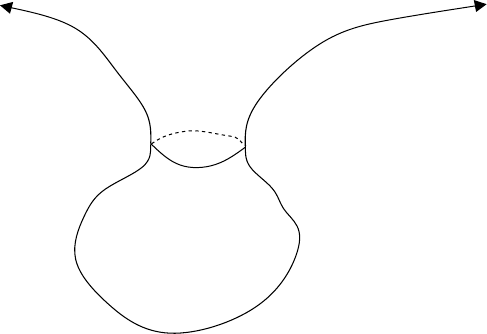}
\put(-130,90){$\infty$}
\put(-90,80){$\tiny{M^3_{\Sigma(t)}}$}
\put(-80,30){$\tiny{\Sigma(t)}$}
\caption{($M^3,g_t$)}
\label{fig:left}
\end{subfigure}
\hspace{0.1cm}
\begin{subfigure}{0.3\textwidth}
\centering
\includegraphics[width = \textwidth]{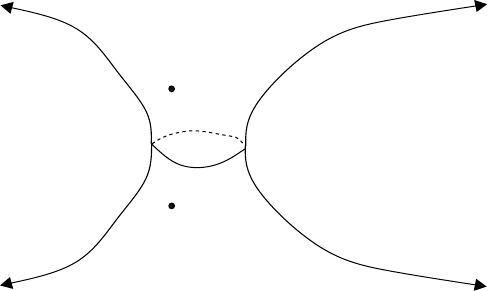}
\put(-86,55){$x$}
\put(-90,13){$\footnotesize{\sigma(x)}$}
\put(-8,-7){$-\infty$}
\put(-5,80){$\infty$}
\caption{$(\bar{M}^3_{\Sigma(t)},\bar{g}_t$)}
\label{fig:middle}
\end{subfigure}
\begin{subfigure}{0.3\textwidth}
\centering
\includegraphics[width = \textwidth]{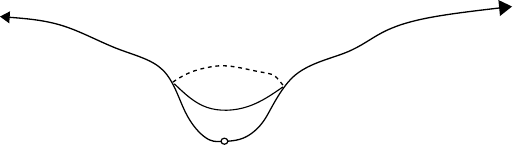}
\put(-10,40){$\infty$}
\put(-80,-8){$-\infty$}
\vspace{0.2cm}
\caption{$(\widetilde{M}^3_{\Sigma(t)},\widetilde{g}_t$)}
\label{fig:right}
\end{subfigure}
\caption{}
\label{fig:combined}
\end{figure}

For a clear visualization of the conformal flow, refer to Figures 1 and 2 in \cite{bray01}.    
We use the notations $|\cdot|$, $\nabla$, $\Delta$ for norm, connection and Laplacian with respect to the original metric $g$, while $\bar{\nabla}$, $\widetilde{\nabla}$, $\bar{\Delta}$ and
$\widetilde{\Delta}$  are the connections and Laplacians with respect to $\bar{g}_t$ and $\widetilde{g}_t$. Here are some properties about the conformal flow from \cite[Theorem 2 and 4]{bray01}, also see \cite[Theorem 2.2]{bray2009riemannian}. 

\begin{theorem} \label{property}
    \begin{enumerate}
    \item For $t_2>t_1\ge 0$, $\Sigma(t_2)$ encloses $\Sigma(t_1)$ without touching it.
        \item The area of $\Sigma(t)$ does not change under the conformal flow:  $|\Sigma(t)|_{g_t}=|\Sigma(0)|_{g_0}$.
        \item If the scalar curvature is nonnegative, the mass of $(M^3,g_t)$ is decreasing.
        \item If $(M^3,g_0)$ is harmonically flat, then for sufficiently large $t$, there exists a diffeomorphism $\phi_t$ between $(M^3,g_t,\Sigma(t))$ and a fixed Schwarzschild metric $(\R^3\setminus\{0\},s)$ such that $|g_t-\phi_t^*(s)|\to 0$ and $|m(t)-m_s|\to 0$,   when $t\to \infty$.
    \end{enumerate}
\end{theorem}
\begin{proof}
    Statement 1 and 2  follows from \cite[Theorem 2 and 3]{bray01}. 

3. We give a short scratch of the proof of Statement 3 for further purpose, and the detail proof is given in \cite[Theorem 3]{bray01}.   Let $m(t)$ be the ADM mass with respect to $(M^3,g_t)$.
Then $-\frac{1}{2}m'(t)$ can be interpered as the ADM mass of an asymptotically flat manifold $(\widetilde{M}^3_{\Sigma(t)},\widetilde{g}_t)$: $\widetilde{m}(t)$ (also in \cite[Lemma 2.7]{bray2009riemannian}), where  
$(\widetilde{M}^3_{\Sigma(t)},\widetilde{g}_t)$ is obtained by the following procedures. 
We glue two copies of $M^3_{\Sigma(t)}$ to get an asymptotically flat manifold  $(\bar{M}^3_{\Sigma(t)},\bar{g}_t)$ with two ends, see the Figure \ref{fig:middle}. The mean curvature matches on both sides of $\Sigma(t)$ as it is a minimal surface. Then we conformally deform the metric to close an end with a harmonic conformal factor $\bar{\phi}$ to obtain $(\widetilde{M}^3_{\Sigma(t)},\widetilde{g}_t)$, 
where $\widetilde{g}_t=\bar{\phi}_t^4\bar{g}_t$ and $\bar{\phi}_t$ is a harmonic function on $(\bar{M}^3_{\Sigma(t)},\bar{g}_t)$ defined in \cite[Equation (76)]{bray01}:
\begin{equation}\label{bphit}
\Delta_{\bar{g}_t}\bar{\phi}_t =0,
        \quad        \lim_{x\to \infty}\bar{\phi}_t=1 \quad\text{and}\quad
        \quad  \lim_{x\to -\infty}\bar{\phi}_t=0,  
\end{equation}
here $\infty$ is the chosen end and $-\infty$ is the other end of $(\bar{M}^3_{\Sigma(t)},\bar{g}_t)$. By symmetry, $\bar{\phi}_t=\frac{1}{2}$ on $\Sigma_t$.  Since the original metric is harmonically flat, the $-\infty$ end in $(\bar{M}^3_{\Sigma(t)},\bar{g}_t)$ becomes a removable singularity in $(\widetilde{M}^3_{\Sigma(t)},\widetilde{g}_t)$, see Figure \ref{fig:right}. 

Let $\phi_t=u_t\bar{\phi}_t$ and let $\widetilde{R}_t$ be the scalar curvature of $(\widetilde{M}^3_{\Sigma(t)},\widetilde{g}_t)$. Thus, denote $\phi_t=u_t\bar{\phi}_t$ and $\phi_t$ satisfies
\begin{equation}
    \Delta \phi_t=0;\quad \lim_{x\to \infty}\phi_t=e^{-t};\quad \phi_t=\frac{u_t}{2} \textit{ on } \Sigma_t.  
\end{equation}
If the scalar curvature of $(M^3,g_0)$ is nonnegative, then the scalar curvature of $(\tilde{M}^3_{\Sigma_t},\tilde{g}_t)$ is also nonnegative, which is due to the formula below,
\begin{equation}
    \widetilde{R}_t=\phi^{-4}_t R_0. 
\end{equation}
Although the manifold is only Lipschitz along the surface $\Sigma(t)$, we apply the positive mass theorem with corners \cite{Pengzi,SvenPengziYTin} to prove that $\tilde{m}_t$ is positive. Hence, the mass of $(M^3,g_t)$ is decreasing.

4.   In this statement, we do not assume nonnegative scalar curvature. However, since $(M^3,g_0)$ is harmonically flat, there exists a compact set $K$ such that the scalar curvature is zero outside $K$. Note that for sufficiently large $t>0$, $\Sigma(t)$ would enclose $K$ which was proven in \cite[Theorem13]{bray01}. In the proof of Theorem 13, the only place using the nonnegativity of scalar curvature is 
Theorem 11.  We can modify the statement and proof of Theorem 11 as below.  

Let $\Sigma_i(t)$ be a connected componment of $\Sigma(t)$, then
$\Sigma(t)$ is a strictly outer minimizing  surface in $(M^3,g_t)$ implies the second variation of the area is nonnegative. Hence, it follows
\begin{equation}
    \int_{\Sigma_i(t)}-\frac{1}{2}R+\frac{1}{2}R^\Sigma-\frac{1}{4}(\lambda_1-\lambda_2)^2\ge0,
\end{equation}
where $\frac{1}{2}R^\Sigma$ is the Gauss curvature of $\Sigma_i(t)$, $\lambda_1$ and $\lambda_2$ are the principal curvature of $\Sigma_i(t)$. Therefore, using Gauss-Bonnet formula and $|\Sigma_i(t)|_{g_t}\le |\Sigma(t)|_{g_t}=A$, we have
\begin{equation}
    \int_{\Sigma_i(t)}(\lambda_1-\lambda_2)^2 \le 16\pi-\int_{\Sigma_i(t)}2R\le   16\pi+2A\max_{K}{|R|}<\infty,
\end{equation}
then the proof of Theorem 13 still works without assuming nonnegative scalar curvature. 
Hence, Theorem 4 in \cite{bray01} implies statement 4.  
\end{proof}

\section{The Harmonic Level Set Formula}
The harmonic level set method, developed by Stern in \cite{Stern}, was subsequently utilized to establish the Positive Mass Theorem by Bray, Kazaras, Khuri, and Stern in \cite{BKKS}. Hirsch, Miao, and Tsang later extended the mass formula to include manifolds with corners in \cite{SvenPengziYTin}. The Positive Mass Theorem with corners was initially proven in \cite{Pengzi}. The following theorem provides a lower bound for the mass in dimension $3$ using a harmonic function which is a weaker version of \cite[Theorem 1.2]{SvenPengziYTin}.
\begin{theorem}\label{ls}(\cite{BKKS, SvenPengziYTin})
      Let $(M^3_{out},g)$  be a complete asymptotically flat manifold with boundary $\Sigma$. Suppose $(\Omega, g_{\Omega})$ is a fill-in of $\Sigma$ such that $g_{\Omega}|_{\Sigma}=g|_{\Sigma}$ and the mean curvatures, with respect to the normal vector pointing outward of $\Omega$, equal from both sides. Let $M^3=M^3_{out}\cup \Omega$. 
      If $H_2(M^3)=0$, then 
\begin{equation} \label{levelset}
            m\ge \frac{1}{16\pi}\int_{M^3}\frac{|\nabla^2 U|^2_{g}}{|\nabla U|_{g}}+R_g|\nabla U|_gdV_g,
    \end{equation}
        where  $U$ is a harmonic function which is asymtotic to one of the asymtotically flat coordinate,  $m$ is the ADM mass of $(M^3,g)$ and $R_g$ is the scalar curvature of $(M^3,g)$.
\end{theorem}
Now we apply Equation \eqref{levelset} to $(\widetilde{M}^3_{\Sigma_t},\widetilde{g}_t)$. To satisfy the topological assumption $H_2(M^3)=0$ in Theorem \ref{ls}, we assume the manifold $(M^3,g)$ only has one single connected outermost minimal surface. Therefore, the outermost minimal surface of $(M^3_{\Sigma_t},g_t)$ is connected, then $\widetilde{M}^3_{\Sigma_t}$ has trivial topology.
\begin{proof}[Proof of Theorem \ref{harmonic theorem}]
We use the notations in Theorem \ref{property}. 
    According to  Theorem \ref{property}, 
    \begin{equation}
    m'(t)=-2\widetilde{m}(t)\quad \text{ and } \quad \lim_{t\to\infty} m(t)=\sqrt{\frac{A}{16\pi}},    
    \end{equation}
    then
    \begin{equation} \label{lsf}
        m=m(0)=\lim_{t\to \infty}m(t)-\int_0^tm'(t)dt=\sqrt{\frac{A}{16\pi}}+\int_0^\infty 2\tilde{m}(t)dt. 
    \end{equation}
Since $\phi_t\to e^{-t}$ at $\infty$, $\widetilde{g}_t=\phi_t^4g$ is asymptotic to $e^{-4t}\delta$ at $\infty$. Therefore, without loss of generality, a harmonic function $\widetilde{\mathcal{U}}_t$ on $(\widetilde{M}^3_{\Sigma(t)}, \widetilde{g}_t)$ asymptotic to a linear function with unit length gradient should satisfy
\begin{equation}
    \widetilde{\mathcal{U}}_t\to e^{-2t}x_1 \text{ at }\infty,
\end{equation}
where we choose $x_1$ to be the first  component of the asymptotic coordinate of $x$.
    
   Applying Theorem \ref{ls} to Equation \eqref{lsf}, we have
    \begin{equation} \label{lsm}
        m\ge \sqrt{\frac{A}{16\pi}}+\frac{1}{8\pi}\int^\infty_0\int_{\widetilde{M}^3_{\Sigma_t}}
    \frac{|\widetilde{\nabla}^2 \widetilde{\mathcal{U}}_t|^2_{\widetilde{g}_t}}{|\widetilde{\nabla} \widetilde{\mathcal{U}}_t|_{\widetilde{g}_t}}+\widetilde{R}_t|\widetilde{\nabla} \widetilde{\mathcal{U}}_t|_{\widetilde{g}_t}dV_{\widetilde{g}_t}dt,
    \end{equation}
    where   $\widetilde{\nabla}$ and $\widetilde{R}_t$ are the connection and the scalar curvature on $(\widetilde{M}^3_{\Sigma(t)},\widetilde{g}_t)$.

    We need to express the inequality \eqref{lsm} under the original metric $g$. Note that on $(\bar{M}^3_{\Sigma(t)},\bar{g}_t)$, due to Lemma \ref{harmclass},   we have $ \overline{\mathcal{U}}_t:=\bar{\phi}_t\widetilde{\mathcal{U}}_t$ satisfying $\Delta_{\bar{g}_t} \overline{\mathcal{U}}_t=0$. Using the asymptotic of $\bar{\phi}_t$ given in Equation \eqref{bphit}, the asymptotics of  $\overline{\mathcal{U}}_t$ are $\overline{\mathcal{U}}_t\to e^{-2t}x_1$ at $\infty$ and $\overline{\mathcal{U}}_t\to 0$ at $-\infty$. 

    Let $\sigma_t$ be the reflection map along $\Sigma(t)$ on $(\bar{M}^3_{\Sigma(t)},\bar{g}_t)$, then we define two harmonic functions:
    \begin{equation} \label{pqt}
        \bar{p}_t(x):=\overline{\mathcal{U}}_t(x)+ \overline{\mathcal{U}}_t(\sigma_t(x)) \quad\text{and}\quad
         \bar{q}_t(x):=\overline{\mathcal{U}}_t(x)- \overline{\mathcal{U}}_t(\sigma_t(x)).
    \end{equation}
     Hence, $\overline{\mathcal{U}}_t=\frac{1}{2} (\bar{p}_t+\bar{q}_t)$,  $\bar{p}_t(x)=\bar{p}_t(\sigma_t(x))$ and $\bar{q}_t(x)=-\bar{q}_t(\sigma_t(x))$. Moreover, the asymptotics of $\bar{p}_t$ and $\bar{q}_t$ are
\begin{align}
\bar{p}_t \to e^{-2t}x_1 \text{ at } \pm\infty;
\quad \bar{q}_t \to\pm e^{-2t}x_1 \text{ at } \pm\infty.
\end{align}
 On $\Sigma_t$, let $\bar{\nu}_t$ is a unit normal vector pointing towards to the chosen end with respect to $\bar{g}_t$. By the definition of $\bar{p}_t$ and $\bar{q}_t$ in \eqref{pqt}, we have
 \begin{equation}
     \bar{\nu}_t(\bar{p}_t)=0 \quad \text{and} \quad \bar{q}_t=0.
 \end{equation}

On $(M^3_{\Sigma(t)},g)$, let $p_t:=u_t\bar{p}_t$ and $q_t:=u_t\bar{q}_t$ be the harmonic functions and  let $\nu_t:=u^2_t\bar{\nu}_{t}$ be a unit vector, we have 
 \begin{equation}
     \nu_t(p_t)=u^2_t\bar{\nu}_t(u_t\bar{p}_t)=u_t^2\bar{p}_t\bar{\nu}_t(u_t)=p_tu_t^{-1}\nu_t(u_t).
 \end{equation}
 Hence, from the definition of $p_t$, $q_t$ and Lemma \ref{harmclass}, it follows that
 \begin{align}
 \begin{split}
     \Delta p_t=0,&\quad p_t=e^{-3t}x_1 \text{ at } \infty \text{ and }  \nu_t(p_t)=p_t\nu_t(\log u_t) \text{ at }\Sigma_t;
     \\ \Delta q_t=0,&\quad q_t=e^{-3t}x_1 \text{ at } \infty \text{ and }  q_t=0 \text{ at }\Sigma_t.
 \end{split}
 \end{align} 
Furthermore, we can use $\phi_t$, $p_t$ and $q_t$ to express $\widetilde{\mathcal{U}}_t$:
\begin{equation} \label{Ut}
\begin{split}
    &\widetilde{\mathcal{U}}_t=\bar{\phi}_t^{-1}\overline{\mathcal{U}}_t=\frac{1}{2}\bar{\phi}_t^{-1}(\bar{p}_t+\bar{q}_t)=\frac{1}{2}\phi_t^{-1}(p_t+q_t)=\frac{1}{2}\mathcal{U}_{t,1}, \text{ 
 on  } M^3_{\Sigma(t)}; 
 \\ &\widetilde{\mathcal{U}}_t=\frac{1}{2}\bar{\phi}_t^{-1}(\bar{p}_t+\bar{q}_t)=\frac{1}{2}\phi_t^{-1}(p_t-q_t)\circ \sigma_t=\frac{1}{2}\mathcal{U}_{t,2}\circ\sigma_t, \text{ 
 on  } \widetilde{M}^3_{\Sigma(t)}\setminus M^3_{\Sigma(t)}, 
\end{split}   
\end{equation}
where $\mathcal{U}_{t,1}$ and $\mathcal{U}_{t,2}$ are defined in this theorem.

Note that the formula for the Hessian under the conformal change $g\to\tilde{g}_t=\phi_t^4 g$ is
\begin{equation}
    \widetilde{\nabla}^2_{ij}f= \partial^2_{ij} f-\widetilde{\Gamma}^l_{ij}\partial_l f
    =\nabla^2_{ij} f-2\phi_t^{-1}\partial_j \phi_t\partial_i f-2\phi_t^{-1}\partial_i \phi_t\partial_j f+2\phi_t^{-1}g^{kl}\partial_k \phi_t\partial_l f g_{ij}.
\end{equation}
Then 
\begin{equation} \label{Hess conf}
    |\widetilde{\nabla}^2 \widetilde{\mathcal{U}}_t|^2_{\widetilde{g}_t}=\phi_t^{-8}\left|\nabla^2 \widetilde{\mathcal{U}}_t-2\phi_t^{-1}\nabla \phi_t\otimes\nabla \widetilde{\mathcal{U}}_t-2\phi_t^{-1}\nabla \widetilde{\mathcal{U}}_t\otimes\nabla \phi_t+2\phi_t^{-1}\langle\nabla\phi_t,\nabla \widetilde{\mathcal{U}}_t\rangle g\right|^2.
\end{equation}
Therefore, using formulae under the conformal metric $\widetilde{R}_t=\bar{\phi}^{-4}_t\bar{R}_t=\phi_t^{-4}R$, $dV_{\widetilde{g}_t}=\phi_t^{6}dV$ and $|\widetilde{\nabla} \widetilde{\mathcal{U}}_t |_{\widetilde{g}_t}=\phi_t^{-2}|\nabla \widetilde{\mathcal{U}}_t |$, then
\begin{equation}
   \begin{split}
       \int_{\widetilde{M}^3_{\Sigma(t)}}\widetilde{R}_t|\widetilde{\nabla} \widetilde{\mathcal{U}}_t|_{\widetilde{g}_t}dV_{\widetilde{g}_t}
       =&\int_{\bar{M}^3_{\Sigma(t)}}\bar{R}_t|\bar{\nabla}\widetilde{\mathcal{U}}_t|_{\bar{g}_t}dV_{\bar{g}_t}
       \\=& \int_{M^3_{\Sigma(t)}} \bar{R}_t|\bar{\nabla}\widetilde{\mathcal{U}}_t|_{\bar{g}_t}+ \bar{R}_t|\bar{\nabla}\widetilde{\mathcal{U}}_t(\sigma_t(x))|_{\bar{g}_t}dV_{\bar{g}_t}
       \\=&\frac{1}{2}\int_{M^3_{\Sigma(t)}}R(|\nabla \mathcal{U}_{t,1}|+|\nabla \mathcal{U}_{t,2}|)dV.
   \end{split}
\end{equation}
Combining \eqref{lsf}, \eqref{Ut} and \eqref{Hess conf}, we obtain 
\begin{equation} \label{mutj}
        \small{
    \begin{split}
             m\ge& \sqrt{\frac{A}{16\pi}}+\frac{1}{8\pi}\int_0^\infty  \int_{\widetilde{M}^3_{\Sigma(t)}}\widetilde{R}_t |\nabla \widetilde{\mathcal{U}}_t|dV_{\widetilde{g}_t}dt+\frac{1}{8\pi} \int_0^\infty \int_{\bar{M}^3_{\Sigma(t)}}\frac{|\nabla^2 \widetilde{\mathcal{U}}_t|^2_{\widetilde{g}_t}}{|\nabla \widetilde{\mathcal{U}}_t|_{\widetilde{g}_t}} dV_{\widetilde{g}_t}dt
             \\ =& \sqrt{\frac{A}{16\pi}}+\int_0^\infty\frac{1}{16\pi}\left( \int_{M_{\Sigma(t)}}  R(|\nabla\mathcal{U}_{t,1}|+ |\nabla \mathcal{U}_{t,2}|)dV\right) dt+
             \\ \frac{1}{16\pi} &\sum_{j=1}^2\int_0^\infty \int_{M^3_{\Sigma(t)}}\frac{\left|\nabla^2 \mathcal{U}_{t,j}-2\phi_t^{-1}\nabla \phi_t\otimes\nabla \mathcal{U}_{t,j}-2\phi_t^{-1}\nabla \mathcal{U}_{t,j}\otimes\nabla \phi_t+2\phi_t^{-1}\langle\nabla\phi_t,\nabla \mathcal{U}_{t,j}\rangle g\right|^2}{|\nabla \mathcal{U}_{t,j}|}dVdt.
        \end{split}
    }    \end{equation}
        Since $R\in L^1(M^3)$ and $|\nabla \mathcal{U}_{t,i}|$ are bounded,  we can switch the order of integration of the second term in the second line of \eqref{mutj} which leads to the final inequality.  
        \end{proof}
        In the following example, we compute $Q$ in  Schwarzschild metric.
        \begin{example} \label{CFS}
           Let $g_{ij}=(1+\frac{m}{2r})^4\delta_{ij}$, then $(\widetilde{M}^3_{\Sigma(t)},\widetilde{g}_t)$ is flat and 
           $\widetilde{g}_t=e^{-4t}\delta$. Since $\widetilde{\mathcal{U}}_t$ is the harmonic function under the metric $\widetilde{g}_t$ asymptotic to $e^{-2t}x_1$, we have  
$\widetilde{\mathcal{U}}_t=e^{-2t}x_1$ and $|\widetilde{\nabla}\widetilde{\mathcal{U}}_t|_{\widetilde{g}_t}=1$. 
        According to \cite[Equation (248)]{bray01}, $u_t=(1+\frac{m}{2r})^{-1}(e^{-t}+\frac{m}{2r}e^t)$, 
        when $r\ge\frac{m}{2}e^{2t}$, then $t(x)=\frac{1}{2}\log \frac{2r}{m}$. 
        In the rigidity case, the conformal flow is a reparametrization of the Schwarzschild metric outside the minimal surface. 

   Note that $\bar{\phi}_t$ is the harmonic function on $(\bar{M}^3_{\Sigma(t)}, u_t^4 g)$ and $\bar{\phi}_t\to 1$ at $\infty$, then using Lemma \ref{harmclass},      
   we  have $\bar{\phi}_t$ is a multiple of the inverse of the conformal factor $(1+\frac{m}{2r})u_t$, therefore,  
   $\bar{\phi}_t=(1+\frac{m}{2r}e^{2t})^{-1}$ 
   and $\phi_t=u_t\bar{\phi}_t= (1+\frac{m}{2r})^{-1}e^{-t}$. Since $\bar{g}_t=u^4_tg=(e^{-t}+\frac{m}{2r}e^t)^4\delta$, it follows $\sigma_t:x\to (\frac{m}{2}e^{2t})^2r^{-2}x$, 
        then $\mathcal{U}_{t,1}=2\widetilde{\mathcal{U}}_t=2e^{-2t}x_1$ implies 
   \begin{equation}
       \mathcal{U}_{t,2}=\mathcal{U}_{t,1}\circ \sigma_t=2(\frac{m}{2}e^{2t})^2r^{-2}e^{-2t}x_1=\frac{m^2}{2}e^{2t}r^{-2}x_1.
   \end{equation}  
        Hence, let $\mathring{\nabla}$ be the connection under $\delta$, we have
        \begin{equation}
        \begin{split}
             Q=&\int_0^{\frac{1}{2}\log \frac{2r}{m}} \sum_{j=1}^2|\nabla \mathcal{U}_{t,j}|dt
             \\=& \int_0^{\frac{1}{2}\log \frac{2r}{m}} \left(1+\frac{m}{2r}\right)^{-2}\sum_{j=1}^2|\mathring{\nabla} \mathcal{U}_{t,j}|dt
             \\=& \int_0^{\frac{1}{2}\log \frac{2r}{m}} \left(1+\frac{m}{2r}\right)^{-2}\left(2e^{-2t}+\frac{m^2}{2}e^{2t}r^{-2}\right)dt
             \\=& \left(1+\frac{m}{2r}\right)^{-1}\left(1-\frac{m}{2r}\right).
        \end{split}
        \end{equation}
        We obtain the same answer for $Q(x)$ as Example \ref{imcf}.
        \end{example}
\section{Mass formula using spinors}
 Here is the mass formula with spinors in the Riemannian case. 
\begin{theorem}(\cite{Witten81, PT81}) Let $M^n$ be a spin asymptotically flat manifold with nonnegative scalar curvature $R$ and total mass $m$ for a chosen end. Then there exists a harmonic spinor $\psi$ satisfying $\mathcal{D}\psi=0$ and  asymtotic to a unit constant spinor at the chosen end. Moreover,
        \begin{equation}
        \label{spin}
            m=\frac{1}{16\pi}\int_{M^n}4|\nabla\psi| ^2+R |\psi|^2 dV.
        \end{equation}
\end{theorem}

The conformal change of metric will induce an isometry between spinor bundles. Since we have two types of conformal changes: 
\[\Phi_t:(M^3_{\Sigma(t)},g)\to (\bar{M}^3_{\Sigma(t)},
\bar{g}_t=u_t^4g)\quad and \quad \bar{\Phi}_t:(\bar{M}^3_{\Sigma(t)},\bar{g}_t)\to(\widetilde{M}^3_{\Sigma(t)},\widetilde{g}_t=\bar{\phi}^4_t\bar{g}_t),\] 
$\Phi_t$, $\bar{\Phi}_t$ would induce two maps between the spinor bundles, and we 
 still denote as  $\Phi_t$, $\bar{\Phi}_t$: 
\begin{equation}
    \Phi_t:(M^3_{\Sigma(t)}, \mathcal{S})\to (\bar{M}^3_{\Sigma(t)}, 
    \mathcal{S})\quad and \quad   \bar{\Phi}_t:(\bar{M}^3_{\Sigma(t)}, \mathcal{S})\to (\widetilde{M}^3_{\Sigma(t)}, \mathcal{S}).
\end{equation}
These maps $\Phi_t$,  $\bar{\Phi}_t$ between spinor bundles are almost identical maps, while the actions of the vector fields on the spinor bundles are different. Let $X$ be a tangent vector on $M^3_{\Sigma_t}$, then for $\psi\in(M^3_{\Sigma(t)}, \mathcal{S})$ and $\bar{\psi}\in(\bar{M}^3_{\Sigma(t)}, \mathcal{S})$, we have
\begin{equation}
    \Phi_t(X\cdot \psi)=u_t^{-2}X\cdot \Phi_t( \psi)\quad \text{and}\quad \bar{\Phi}_t(X\cdot \bar{\psi})=\bar{\phi}_t^{-2}X\cdot \bar{\Phi}_t(\bar{\psi}).
\end{equation}

Before we prove the main theorem in this section, let us first identify the action induced by the reflection on the spinor bundle of $\Sigma(t)$. 
\begin{lemma} \label{pm iv}
    Let $\bar{M}^3_{\Sigma}$ be the asymptotically flat manifold with two ends with an outermost minimal surface ${\Sigma}$. Suppose there exists a reflection map $\sigma$ along $\Sigma$, i.e.,  $\bar{M}^3_{\Sigma}$ is symmetric along $\Sigma$ (see Figure \ref{fig:middle}).  Then for a spinor $\psi$ on $\bar{M}^3_{\Sigma}$,
    \begin{equation}
        \sigma^*\psi|_{\Sigma}=i\nu\cdot \psi|_{\Sigma} \text{ or } -i\nu\cdot \psi|_{\Sigma},
    \end{equation}
    where $\nu$ is the unit normal vector on $\Sigma$ and $\sigma^*$ is the isomorphism of the spinor bundle on $\Sigma$  induced by $ \sigma$.
\end{lemma}
\begin{proof} To restrict a spinor on a hypersurface, we follow the presentation in \cite[Page 6]{harmonic spinor}.  
 The Clifford algebra $\operatorname{Cl}(\R^n)$ can be splited into odd and even part: $\operatorname{Cl}(\R^n)=\operatorname{Cl}^0(\R^n)\oplus \operatorname{Cl}^1(\R^n)$. Then there is an algebra isomorphism $\operatorname{Cl}(\R^{n-1})\to \operatorname{Cl}^0(\R^{n})$ induced by $e_i\to e_i\cdot e_n$, where $\{e_1,\dots,e_n \}$ is an orthonormal frame on $\R^n$.  

In this lemma, $n=3$ is odd. Let $\{e_1,e_2,e_3\}$ be a orthonormal frame on $M^3$. In particular, $e_3=\nu$ on $\Sigma$.  Let $\mathcal{S}_M$ and $\mathcal{S}_\Sigma$ be the spinor bundle on $\bar{M}^3_{\Sigma}$ and $\Sigma$, then there is a isomorphism between the spinor bundles $\mathcal{S}_M|_{\Sigma}$ and $\mathcal{S}_\Sigma$. Since $i\nu$ is an automorphism on $\mathcal{S}_\Sigma$ and $(i\nu)^2=1$, we may decompose $\mathcal{S}_\Sigma$ into the eigenspace $\mathcal{S}^+_\Sigma$ and  $\mathcal{S}^-_\Sigma$ with respect to the eigenvalue $1$ and $-1$. As $i\nu$ anticommutes with $e_1\nu$, $e_2\nu$ and commutes with $e_1e_2$, then the action of $i\nu$ on $\mathcal{S}_\Sigma$ is $ie_1e_2$ or $-ie_1e_2$. 

Next, we study the action of $ie_1e_2$ on the spinor bundle of $\Sigma$.  
Since on $\mathcal{S}_\Sigma$, the complex spinor bundle can be identified as 
$\operatorname{Span}_{\mathbb{C}}\{\mathbf{1},\eta:=\frac{1}{\sqrt{2}}(e_1e_3-ie_2e_3)\}$, see \cite[Page 79]{Jost},  the actions of $e_1e_3$ and $e_2e_3$ are
\begin{equation}
    (e_1e_3)\cdot\mathbf{1}=\eta,\quad 
    (e_1e_3)\cdot\mathbf{\eta}=-\mathbf{1},\quad 
    (e_2e_3)\cdot\mathbf{1}=i\eta,\quad 
    (e_2e_3)\cdot\mathbf{\eta}=i\mathbf{1}.
\end{equation}
Note that $ie_1e_2=ie_1e_3e_2e_3$, then
\begin{equation}
    (ie_1e_2)\cdot \mathbf{1}=\mathbf{1} \quad \text{and} \quad 
    (ie_1e_2)\cdot\eta=-\eta
\end{equation}
Finally, we have $\sigma^*(e_1)=e_1$, $\sigma^*(e_2)=e_2$ and $\sigma^*(e_3)=-e_3$, then $\sigma^*(\mathbf{1})=\mathbf{1}$ and $\sigma^*(\eta)=-\eta$. Hence,  
the action of $\sigma^*$ on $\mathcal{S}_\Sigma$ is the same as $ie_1e_2$.
\end{proof}
\begin{remark}
    The $\pm$ sign depends on the choice of the spin group action on the spinor bundle.
\end{remark}

\begin{proof}[Proof of Theorem \ref{spinor theorem}]
    We use Equation \eqref{spin} to express $\widetilde{m}(t)$, then we apply Equation \eqref{lsf} to achieve a type of inequality \eqref{1}.  Let $\widetilde{\psi}_t$ be the harmonic spinor on $(\widetilde{M}^3_{\Sigma(t)},\widetilde{g}_t)$ which asymptotic to a unit constant spinor at $\infty$. We need to rewrite the integral formula of $\widetilde{\psi}_t$ under the original metric $(M^3,g)$. 
    
For simplicity, we denote $\bar{\mathcal{D}}$ and $\widetilde{\mathcal{D}}$ to be the Dirac operator on $(\bar{M}^3_{\Sigma(t)},\bar{g}_t)$ and $(\widetilde{M}^3_{\Sigma(t)},\widetilde{g}_t)$. Since $\widetilde{g}_t=\bar{\phi}^4_tg$, we have the conformal formula for Dirac operator
    \begin{equation}
       \bar{\Phi}_t \bar{\mathcal{D}}(\bar{\phi}_t^{2} \Phi_t^{-1}(\widetilde{\psi}_t))=\bar{\phi}_t^{4}\widetilde{\mathcal{D}}(\tilde{\psi}_t).
    \end{equation}
    Let $\bar{\psi}_t=\bar{\phi}_t^2\bar{\Phi}_t^{-1}(\widetilde{\psi}_t)$
and $\psi_t=u_t^2\Phi_t^{-1}(\bar{\psi}_t)$, then $\bar{\psi}_t$ and $\psi_t$ are harmonic spinors with respect to $\bar{g}_t$ and $g_t$. Moreover, because of the asymptotic of $\phi_t$ in Equation \eqref{bphit}, it follows that the asymptotics of $\bar{\psi}_t$ and $\psi_t$ are
\begin{equation}
    \lim_{x\to \infty}\psi_t=e^{-2t}\psi_0 ; \quad \lim_{x\to \infty}\bar{\psi}_t=\bar{\psi}_0\quad and\quad \lim_{x\to -\infty}\bar{\psi}_t=\mathbf{0},
\end{equation}
where $\psi_0$,  $\bar{\psi}_0$ are unit constant spinors and $\mathbf{0}$ is the zero spinor. Since $\psi_t$ is only defined on $(M^3_{\Sigma(t)},g)$, we need to analyze the boundary behavior of $\psi_t$ on $\Sigma(t)$ below. 

Let $\sigma_t$ be the reflection map $\bar{M}^3_{\Sigma(t)}\to \bar{M}^3_{\Sigma(t)}$, then $\sigma_t$ induces a map $\sigma^*_t$ on the spinor bundle on $\bar{M}^3_{\Sigma(t)}$.  
 Let $\bar{\mathbf{p}}_t$ and $\bar{\mathbf{q}}_t$ be two harmonic spinors on $(\bar{M}^3_{\Sigma(t)},\bar{g}_t)$ defined as follow
\begin{equation}
    \bar{\mathbf{p}}_t=\bar{\psi}_t+\sigma^*_t \bar{\psi}_t,\quad \bar{\mathbf{q}}_t=\bar{\psi}_t-\sigma^*_t \bar{\psi}_t.
\end{equation}
Hence, 
\begin{equation} \label{bpsipq}
    \bar{\psi}(x)=\frac{1}{2}(\bar{\mathbf{p}}_t+\bar{\mathbf{q}}_t)(x), \text{ if } x\in M^3_{\Sigma(t)}; \quad 
     \bar{\psi}(x)=\frac{1}{2}(\sigma_t^*)^{-1}(\bar{\mathbf{p}}_t-\bar{\mathbf{q}}_t)(\sigma_t(x)), \text{ if } x\in \bar{M}^3_{\Sigma(t)}\setminus M^3_{\Sigma(t)}. 
\end{equation}
Without loss of generality, according to Lemma \ref{pm iv}, we assume the action of $\sigma^*_t$ on the spinor bundle of $\Sigma_t$ is $i\bar{\nu}_t\cdot$, then the restriction of $\bar{\mathbf{p}}_t$ and $\bar{\mathbf{q}}_t$ on $\Sigma$ are
\begin{equation}
  \bar{\mathbf{p}}_t=\bar{\psi}_t+i\bar{\nu}_t\cdot \bar{\psi}_t,\quad \bar{\mathbf{q}}_t=\bar{\psi}_t-i\bar{\nu}_t\cdot \bar{\psi}_t.  
\end{equation}
Hence, on $\Sigma(t)$, we have the Lopatinsky-Shapiro type boundary conditions for $\bar{\mathbf{p}}_t$ and $\bar{\mathbf{q}}_t$
\begin{equation}
    i\bar{\nu}_t\cdot \bar{\mathbf{p}}_t=\bar{\mathbf{p}}_t,\quad i\bar{\nu}_t \cdot \bar{\mathbf{q}}_t=-\bar{\mathbf{q}}_t. 
\end{equation}
Recall that $\bar{g}_t=u_t^4 g$,
we set $\mathbf{p}_t=u_t^2\Phi^{-1}(\bar{\mathbf{p}}_t)$ and 
 $\mathbf{q}_t=u_t^2\Phi^{-1}(\bar{\mathbf{q}}_t)$, according to Remark \ref{con spin}, we have $\mathbf{p}_t$ and $\mathbf{q}_t$ are harmonic spinors on $(M^3_{\Sigma(t)},g)$ with boundary conditions 
 \begin{equation}
    i\nu_t\cdot \mathbf{p}_t=\mathbf{p}_t,\quad i\nu_t \cdot \mathbf{q}_t=-\mathbf{q}_t, \textit{ on } \Sigma(t),
\end{equation}
and 
\begin{equation}
    \lim_{x\to \infty} \mathbf{p}_t=\lim_{x\to \infty} u^2_t\Phi^{-1}(\bar{p}_t)=u^2_t\Phi^{-1}(\bar{\psi}_0)= e^{-2t}\psi_0, \quad \lim_{x\to \infty} \mathbf{q}_t=e^{-2t}\psi_0.
\end{equation}

Now we can use $\mathbf{p}_t$ and $\mathbf{q}_t$ to express $\widetilde{\psi}_t$. Since $\bar{\phi}_t(\sigma(x))=1-\bar{\phi}_t(x)$ and $\phi_t(x)=u_t(x)\bar{\phi}_t(x)$, when $x\in M^3_{\Sigma(t)}$, according to \eqref{bpsipq} and Remark \ref{con spin}, 
\begin{equation} \label{ab}
   \begin{split}
       \widetilde{\psi}_t(x)=&\frac{1}{2}\phi_t^{-2}(\bar{\Phi}_t\circ\Phi_t)(\mathbf{p}_t+\mathbf{q}_t)(x), \text{ if }x\in M^3_{\Sigma(t)};
       \\ \widetilde{\psi}_t(x)=&\frac{1}{2}(\bar{\Phi}_t\circ (\sigma^*_t)^{-1}\circ\Phi_t)[(u_t-\phi_t)^{-2}(\mathbf{p}_t-\mathbf{q}_t)](\sigma_t(x)), \text{ if }x\in \widetilde{M}^3_{\Sigma(t)}\setminus M^3_{\Sigma(t)}.
   \end{split} 
\end{equation}
 
Let $\{e_1,e_2,e_3\}$ be an  orthonormal tangent frame for $(M^3,g)$ and let $\widetilde{e}_j=\phi_t^{-2}e_j$. From Equation \eqref{gt},
\begin{equation} \label{gt1}
   |\widetilde{\nabla}\widetilde{\psi}_t|_{\widetilde{g}_t}=\sum_{j=1}^3\phi^{-2}_t\left|\nabla_{e_j}(\phi_t^{-2}\psi_t)-[\nabla_{e_j} \log (\phi_t)+{e_j}\cdot \nabla \log \phi_t](\phi_t^{-2}\psi_t)\right| \text{ on } M^3_{\Sigma(t)}. 
\end{equation}

On $\widetilde{M}^3_{\Sigma(t)}\setminus M^3_{\Sigma(t)}$, for simplicity, let $\phi_{t,1}:=u_t-\phi_t$, $\psi_t^1:=\frac{1}{2}(\mathbf{p}_t-\mathbf{q}_t)$, then combining the second equation in \eqref{ab} and \eqref{gt1},  we have  
\begin{equation} \label{gt2}
\begin{split} 
    &|\widetilde{\psi}_t|_{\widetilde{g}_t}
    =\phi_{t,1}^{-2}|\psi_t^1|,
\\ &|\widetilde{\nabla}\widetilde{\psi}_t|_{\widetilde{g}_t}(x)=\sum_{j=1}^3\phi^{-2}_{t,1}\left|\nabla_{e_j}(\phi_{t,1}^{-2}\psi_t^1)-[\nabla_{e_j} \log (\phi_{t,1})+{e_j}\cdot \nabla \log \phi_{t,1}](\phi_{t,1}^{-2}\psi_{t}^1)\right|.
\end{split}
\end{equation}

Let $\psi_t^0:=\frac{1}{2}(\mathbf{p}_t+\mathbf{q}_t)=\psi_t$ and $\phi_{t,0}:=\phi_t$.
Applying the conformal formulae $dV_{\widetilde{g}_t}=\phi^6_tdV_{g_0}$, $\widetilde{R}_t=\phi_t^{-4}R$, then plugging Equation \eqref{gt1} and \eqref{gt2} in \eqref{spin},  we have
\begin{equation} \label{mt2}
\small{
    \begin{split}
     \widetilde{m}_t=&\frac{1}{16\pi}\int_{\widetilde{M}^3_{\Sigma_t}}4|\nabla\widetilde{\psi}_t|_{\widetilde{g}_t} ^2+\widetilde{R}_t |\widetilde{\psi}_t|^2_{\widetilde{g}_t} dV_{\widetilde{g}_t}
     \\=& \sum_{l=0}^1\frac{1}{16\pi}\int_{M_{\Sigma_t}}4\phi_{t,l}^2\sum_{j=1}^3\left|\nabla_{e_j}(\phi_{t,l}^{-2}\psi_t^l)-[\nabla_{e_j} \log \phi_{t,l}+{e_j}\cdot \nabla \log \phi_{t,l}](\phi_{t,l}^{-2}\psi_t^l)\right|^2
     +R \phi_{t,l}^{-2} |\psi^l_t|^2 dV.
    \end{split}
    }
\end{equation}
Hence, using Equation \eqref{lsf}, we integrate \eqref{mt2} to obtain the  final Equation \eqref{2} with $P(x)$ and $Q(x)$ given in \eqref{PQ2}.
\end{proof}

This paper was written for a special issue in memory of Robert Bartnik (1956-2022). The first author had the pleasure of calling Robert Bartnik his friend and role model. I met Robert the first time when I was 27, fresh out of graduate school. Robert was a legend in the field; it was quite an honor to meet him, which I got to do for two weeks in Australia. Robert went out of his way to help me achieve my goals, however he could. He was very generous with his time and suggested important papers for me to read. He was friendly, hilarious at times, and fun to be around, both mathematically and as a human being. Robert Bartnik's mathematical contributions are very important, but the example he set as a person on how to build community is one of the reasons geometric relativity is thriving today.

\vspace{1in}

\appendix 
\section{Formulae under conformal change of metrics}
We list some formulae for spinors and harmonic functions under conformal change of metrics on an 
 $n$-dimensional 
 manifold $M^n$.

Suppose $g'=e^{2u}g$ on $M^n$. Let $\mathcal{E}=\{e_1,...,e_n\}$ be an orthonormal tangent frame which associate with the orthonormal tangent frame $\{e_1'...,e_n'\}$ and $e'_i=e^{-u} e_i$. The map $\Psi$ can be lifted to a map between principal $\text{Spin}_n$ bundles which we denote $\Phi$. See \cite[Page 133]{Lawson} and \cite[Page 308]{GSP}, $\Phi$ is an isometry. 

\begin{lemma} 
    Suppose $\varphi'$ is a spinor on $(M^n,g')$, then

\begin{equation}
   \mathcal{D}'\varphi'= e^{-\frac{n+1}{2}u}\Phi\left[ \mathcal{D}(e^{\frac{n-1}{2}u}\Phi^{-1}\varphi')\right].
\end{equation}
\end{lemma}
\begin{proof}
Suppose $\varphi'=\Phi(\varphi)$, according to \cite[Lemma 5.27]{Lawson},
\begin{equation}\label{connection}
    \nabla'_V\varphi'=\Phi \left[\nabla_V\varphi-\frac{1}{2}(\nabla_V u+V\cdot \nabla u)\varphi\right].
\end{equation}
Hence, 
\begin{equation}
    \begin{split}
        {\mathcal{D}}'\varphi'=&\sum_{j=1}^n\Phi\left\{e_j\cdot[\nabla_{e'_j}    \varphi-\frac{1}{2}(\nabla_{e'_j} u+{e'_j}\cdot \nabla u)\varphi]\right\}
        \\=& \Phi\left\{e^{-u}[\mathcal{D}\varphi+\frac{n-1}{2}(\nabla w)\cdot\varphi]\right\}
        \\=&e^{-\frac{n+1}{2}u}\Phi\left[\mathcal{D}(e^{\frac{n-1}{2}u}\varphi)\right]
    \end{split}
\end{equation}
    
\end{proof}
\begin{remark} \label{con spin}
 In dimension $3$, suppose $\widetilde{g}_t=\phi_t^4g$, let $\Phi_t(\psi_t)=\phi^2_t\widetilde{\psi}_t$, if $\widetilde{\mathcal{D}}(\widetilde{\psi}_t)=0$, then $\mathcal{D}(\psi_t)=0$.
\end{remark}

Let $\{e_1,e_2,e_3\}$ be a orthonormal tangent frame for $(M^3,g)$. We denote $\widetilde{e}_j=\phi_t^{-2}e_j$. From Equation \eqref{connection},
\begin{equation}
    \widetilde{\nabla}_{\widetilde{e}_j}\widetilde{\psi}_t=\Phi_t\left\{\nabla_{\widetilde{e}_j}(\phi_t^{-2}\psi_t)-\frac{1}{2}[\nabla_{\widetilde{e}_j} 2\log (\phi_t)+{\widetilde{e}_j}\cdot \nabla 2\log (\phi_t)](\phi_t^{-2}\psi_t)\right\},
\end{equation}
then we have the norm of the gradient:
\begin{equation} \label{gt}
   |\widetilde{\nabla}\widetilde{\psi}_t|_{\widetilde{g}_t}=\sum_{j=1}^3\phi^{-2}_t\left|\nabla_{e_j}(\phi_t^{-2}\psi_t)-[\nabla_{e_j} \log (\phi_t)+{e_j}\cdot \nabla \log \phi_t](\phi_t^{-2}\psi_t)\right|. 
\end{equation}

In addition, we also use the property of harmonic functions under conformal metrics.
Here is the lemma from \cite[Page 69]{bray01}.
\begin{lemma} \label{harmclass}
    Let $g_2$ and $g_1$ be two conformal metrics on an $n$-dimensional manifold $M^n$, $n\ge 3$. Suppose $g_2=u^{\frac{4}{n-2}}g_1$, then for any smooth function $\phi$, 
    \begin{equation}
        \Delta_{g_1}(u\phi)=u^{\frac{n+2}{n-2}}\Delta_{g_2}\phi+\phi\Delta_{g_1}u.        
    \end{equation}
    In particular, if $\Delta_{g_1}u=0$ and $\Delta_{g_2}\phi=0$, then $\Delta_{g_1}(u\phi)=0$.
\end{lemma}


\begin{thebibliography}{99}
\bibitem{bray01} Hubert Bray,
\textit{Proof of the Riemannian Penrose inequality using the positive mass theorem}, J.  Diff. Geom.,  59.2 (2001): 177-267.
\bibitem{Lawson}H. Blaine Lawson and Marie L Michelsohn,  Spin Geometry. Princeton Mathematical Series, vol. 38. Princeton
University Press, Princeton (1989)
\bibitem{GSP} Colin Guillarmou, Sergiu Moroianu and
Jinsung Park,
\textit{Bergman and Calderón Projectors for Dirac Operators},  J Geom Anal (2012)

\bibitem{bray2009riemannian} Hubert Bray and Dan Lee
 \textit{On the Riemannian Penrose inequality in dimensions less than $8$}, Duke Math. J. 148 (2009) 81

\bibitem{BKKS} Hubert Bray, Demetre Kazaras,  Marcus Khuri and
    Daniel Stern,
   \textit{Harmonic functions and the mass of 3-dimensional asymptotically
   flat Riemannian manifolds},
   J. Geom. Anal, 2022.

\bibitem{SY79} Richard Schoen and Shing Tung Yau, \textit{On the proof of the positive mass conjecture in general relativity}, Comm. Math. Phys., 65(1):45–76, 1979.
\bibitem{Witten81}  Edward Witten,  \textit{A new proof of the positive energy theorem},  Comm. Math.
Phys., 80(3):381–402, 1981.

\bibitem{Stern} D. Stern, \textit{Scalar curvature and harmonic maps to $S
^1$}, J. Diff. Geom., 2022

\bibitem{HI01}  Gerhard Huisken and Tom Ilmanen, \textit{The inverse mean curvature flow and the
Riemannian Penrose inequality}, J. Diff. Geom., 59(3):353–437, 2001.

\bibitem{PT81} Thomas Parker and Clifford Taubes,
\textit{On Witten's proof of the positive energy theorem}, Comm. Math. Phys. 84(1982), no.2, 223–238.

\bibitem{SvenPengziYTin}
Sven Hirsch, Pengzi Miao and Tin-Yau Tsang, \text{Mass of asymptotically flat 3-manifolds with boundary}, Comm. Anal. Geom., 2021.

\bibitem{harmonic spinor}
    Christian B\"{a}r,
   \textit{Metrics with harmonic spinors},
   Geom. Funct. Anal,
   volume 6,
   1996,
   899--942

\bibitem{HanKhuri}
Qing Han and Marcus Khuri, 
\text{The conformal flow of metrics and the general Penrose inequality}, 
Tsinghua Lectures in Mathematics ALM,
pp. 227-242

\bibitem{Jost}Jürgen Jost, Riemannian geometry and geometric analysis, Vol. 42005. Berlin: Springer, 2008.

\bibitem{charge Penrose} Marcus 
 Khuri,  Gilbert Weinstein, Sumio Yamada, 
\textit{Proof of the Riemannian Penrose inequality with charge for multiple black holes}, J. Diff. Geom.106(2017), no.3, 451–498

\bibitem{Penrose Survey} Marc Mars,
\textit{Present status of the Penrose inequality}, Classical Quantum Gravity, 26(2009), no.19, 193001, 59 pp.

\bibitem{Pengzi} Pengzi Miao, \textit{Positive mass theorem on manifolds admitting corners along a hypersurface}, Adv. Theor.
Math. Phys., 6 (2002). no.6, 1163-1182.

\bibitem{PDE Penrose}
 Hubert Bray and Marcus  Khuri, \textit{P.D.E.'s which imply the Penrose conjecture},
Asian J. Math. 15 (2011), no. 4, 557–610.
\end{thebibliography}
\end{document}